\documentclass[reqno, 12pt]{amsart}

\usepackage{amssymb}
\usepackage{latexsym}
\usepackage{epsfig}
\usepackage{epsf}
\usepackage{float}
\usepackage{a4wide}
\usepackage{hyperref}

\newtheorem{theorem}{Theorem}
\newtheorem{lemma}{Lemma}
\newtheorem{corollary}{Corollary}

\def\R{{\mathbb R}}  
\def\N{{\mathbb N}}  
\def\P{{\mathbb P}}    
\def\Z{{\mathbb Z}}

\def\E{{\mathbb E}}
\def\F{{\mathbb F}}

\def\I{{\mathbb I}}

\def\bA{{\mathbf A}}
\def\bB{{\mathbf B}}
\def\bC{{\mathbf C}}
\def\bD{{\mathbf D}}

\def\bH{{\mathbf H}}
\def\b0{{\mathbf 0}}
\def\be{{\mathbf e}}
\def\bn{{\mathbf n}}
\def\bv{{\mathbf v}}
\def\bx{{\mathbf x}}
\def\by{{\mathbf y}}
\def\bz{{\mathbf z}}

\def\calV{{\mathcal V}}

\def\calG{{\mathcal G}}
\def\calI{{\mathcal I}}
\def\calC{{\mathcal C}}

\def\calD{{\mathcal D}}
\def\calP{{\mathcal P}}
\def\calW{{\mathcal W}}
\def\calF{{\mathcal F}}

\def\reff#1{(\ref{#1})}
\def\proofof #1{{\noindent \emph{Proof of #1}.}}
\def\endproof{$\square$ \vskip 2mm}

\begin{document}

\title[The time constant and critical probabilities]{The time constant and critical probabilities\\ in percolation models}
\author{Leandro P. R. Pimentel}
\address{Institut de Math\'{e}matiques\\
\'{E}cole Polytechinique F\'{e}d\'{e}rale de Lausanne\\ 
CH-1015 Lausanne\\
Switzerland.} 
\email{leandro.pimentel@epfl.ch}



\begin{abstract}
We consider a first-passage percolation (FPP) model on a Delaunay triangulation $\calD$ of the plane. In this model each edge $\be$ of $\calD$ is independently equipped with a nonnegative random variable $\tau_\be$, with distribution function $\F$, which is interpreted as the time it takes to traverse the edge. Vahidi-Asl and Wierman \cite{VW90} have shown that, under a suitable moment condition on $\F$, the minimum time taken to reach a point $\bx$ from the origin $\b0$ is asymptotically $\mu(\F)|\bx|$, where $\mu(\F)$ is a nonnegative finite constant. However the exact value of the time constant $\mu(\F)$ still a fundamental problem in percolation theory. Here we prove that if $\F(0)<1-p_c^*$ then $\mu(\F)>0$, where $p_c^*$ is a critical probability for bond percolation on the dual graph $\calD^*$.
\end{abstract}

\maketitle 

\section{Introduction}
First-passage percolation theory on periodic graphs was presented by Hammersley and Welsh \cite{HW65} to model the spread of a fluid through a porous medium. In this paper we continue a
study of planar first-passage percolation models on random graphs, initiated by Vahidi-Asl and Wierman \cite{VW90}, as follows. Let $\calP$ denote the set of points realized in a two-dimensional
homogeneous Poisson point process with intensity $1$. To each $\bv\in\calP$
corresponds an open polygonal region $\bC_\bv=\bC_\bv(\calP)$, the Voronoi tile at $\bv$, consisting of the set of points of $\R^2$ which are closer to $\bv$ than to any other $\bv'\in\calP$. Given $\bx\in\R^2$ we denote by $\bv_\bx$ the almost surely unique point in $\calP$ such that $\bx\in\bC_{\bv_\bx}$. The collection $\{ \bC_\bv \,:\, \bv\in\calP\}$ is called the Voronoi Tiling of the plane based on $\calP$. 

\begin{figure}[tb]
\begin{center}
\includegraphics[width=0.5\textwidth]{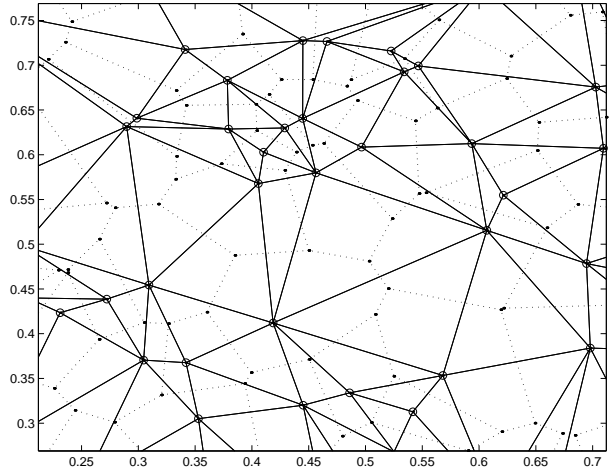}
\end{center}
\caption{The Delaunay Triangulation and the Voronoi Tessellation.}\label{F:vor2}\end{figure}

The Delaunay Triangulation $\calD$ is the graph where the vertex set $\calD_v$ equals $\calP$ and the edge set $\calD_e$ consists of non-oriented pairs $(\bv,\bv')$ such that $\bC_\bv$ and $\bC_{\bv'}$ share a one-dimensional edge (Figure \ref{F:vor2}). One can see that almost surely each Voronoi tile is a convex and bounded polygon, and the graph $\calD$ is a
triangulation of the plane \cite{M91}. The Voronoi Tessellation $\calV$ is the graph where the vertex set $\calV_v$ is the set of vertices of the Voronoi tiles and the edge set $\calV_e$ is the set of edges of the Voronoi tiles. The edges of $\calV$ are segments of the perpendicular bisectors of the edges of $\calD$. This establishes duality of $\calD$ and $\calV$ as planar graphs: $\calV=\calD^*$.

To each edge $\be\in\calD_e$ is
independently assigned a nonnegative random variable $\tau_{\be}$ from a common
distribution $\F$, which is also independent of the Poisson point process that generates $\calP$. From now on we  denote $(\Omega,\calF,\P)$ the probability space induced by the Poisson point process $\calP$ and the passage times $(\tau_\be)_{ \be\in\calD_e }$. The passage time $t(\gamma)$ of a path $\gamma$ in the Delaunay Triangulation is the sum of the passage times of the edges in $\gamma$. The first-passage time between two vertices $\bv$ and $\bv'$ is defined by
\[
T(\bv,\bv'):=\inf\{t(\gamma)\,;\,\gamma\in\calC(\bv,\bv')\}\,,
\]
where $\calC(\bv,\bv')$ the set of all paths connecting $\bv$ to $\bv'$. Given $\bx,\by\in\R^2$ we define $T(\bx,\by):=T(\bv_\bx ,\bv_\by)$.

To state the main result of this work we require some definitions involving a bond percolation model on the Voronoi Tessellation $\calV$. Such a model is constructed by choosing each edge of $\calV$ to be
open independently with probability $p$. An open path is a path composed of open edges. We denote $\P^*_p$ the law induced by the Poisson point process and the random state (open or not) of an edge. Given a planar graph $\calG$ and $\bA,\bB\subseteq\R^2$ we say that a self-avoiding path $\gamma=(\bv_1,...,\bv_k)$ is a path connecting $\bA$ to $\bB$ if $[\bv_1,\bv_2]\cap\bA\neq\emptyset$ and $[\bv_{k-1},\bv_k]\cap\bB\neq\emptyset$ ($[\bx,\by]$ denotes the line segment connecting $\bx$ to $\by$). For $L>0$ let $A_{L}$ be the event that there exists an open path $\gamma=(\bv_j)_{1\leq j\leq h}$ in $\calV$, connecting $\{0\}\times[0,L]$ to $\{3L\}\times[0,L]$, and with $\bv_j\in[0,3L]\times [0,L]$ for all $j=2,\dots,h-1$. In this case we also say that $\gamma$ crosses the rectangle $[0,3L]\times [0,L]$. Define the function 
\[
\eta^*(p):=\liminf_{L\to\infty}\P^*_p(A_L)\,, 
\]
and consider the percolation threshold,
\begin{equation}\label{critical}
p_{c}^*:=\inf\{p>0\,:\,\eta^*(p)=1\}\,.
\end{equation}
We have that $p_c^*\in(0,1)$, which follows by standard arguments in percolation theory. For more in percolation thresholds on Voronoi tilings we refer to \cite{BR-104,BR-204,Z96}.

\begin{theorem}\label{t1}
If $\F(0)<1-p_{c}^*$ then there exist constants $c_j >0$
such that for all $n\geq 1$
\begin{equation}\label{et1}
\P\big(T(\b0,\bn)<c_1 n\big)\leq c_2\exp(-c_3 n)\,,
\end{equation}
where $\b0:=(0,0)$ and $\bn:=(n,0)$.
\end{theorem}

To show the importance of Theorem \ref{t1} we recall two fundamental results proved by Vahidi-Asl and Wierman \cite{VW90,VW92}. Consider the growth process
\[
\bB_\bx (t):=\{\by\in\R^2\,:\,\by\in c(\bC_\bv)\mbox{ with }\bv\in\calD_v\mbox{ and }T(\bv_\bx,\bv)\leq t\}\,.
\] 
where $c(\bC)$ denotes the closure of $\bC\in\R^2$. Set
\[  
\mu(\F):=\inf_{n>0}\frac{\E T(\b0,\bn)}{n}\in[0,\infty]\,.
\]
and let $\tau_1,\tau_2,\tau_3$ be independent random variables with distribution $\F$. If
\begin{equation}\label{c1}
\E\big(\min_{j=1,2,3}\{\tau_j\}\big)<\infty
\end{equation}
then $\mu(\F)<\infty$ and for all unit vectors $\vec{\bx}\in S^1$ ($|\vec{\bx}|=1$) $\P$-a.s.
\begin{equation}\label{e5}
\lim_{n\to\infty}\frac{T(\b0,n\vec{\bx})}{n}=\lim_{n\to\infty}\frac{\E T(\b0,\bn)}{n}=\mu(\F)\,.
\end{equation}
Further, if 
\begin{equation}\label{c2}
\E\big(\min_{j=1,2,3}\{\tau_j\}^2\big)<\infty
\end{equation}
and $\mu(\F)>0$ then for all $\epsilon>0$ $\P$-a.s. there exists $t_0 >0$ such that for all $t>t_0$
\begin{equation}\label{e6}
(1-\epsilon ) t\bD(1/\mu)\subseteq \bB_\b0(t)\subseteq (1+\epsilon )t
\bD(1/\mu)\, ,
\end{equation}
where $\bD(r) :=\{ \bx\in\R^2\, : \, |\bx|\leq r\}$.

We note here that the asymptotic shape is an Euclidean ball due to the statistical invariance of the Poisson point process. Unfortunately the exact value of the time constant $\mu(\F)$, as a functional of $\F$, still a basic problem in first-passage percolation theory. Our result provides a sufficient condition on $\F$ to ensure $\mu(\F)>0$. 
\begin{corollary}\label{col1}
Under assumption \reff{c1}, if $\F(0)<1-p_{c}^*$ then $\mu(\F)\in(0,\infty)$. 
\end{corollary}

\proofof{Corollary \ref{col1}} Together with the Borel-Cantelli Lemma, Theorem \ref{t1} and \reff{e5} imply 
\[ 
0<c_1\leq\liminf_{n\to\infty}\frac{T(\b0,\bn)}{n}=\lim_{n\to\infty}\frac{T(\b0,\bn)}{n}=\mu(\F)\,<\infty \,,
\]
which is the desired result. \endproof

For FPP models on the $\Z^2$ lattice Kesten (1986) have shown that $\F(0)< 1/2=p_c(\Z^2)$ (the critical probability for bond percolation on $\Z^2$) is a
sufficient condition to get \reff{et1} by using a
stronger version of the BK-inequality. Here we follow a different method and
apply a simple renormalization argument to obtain a similar result. We expect that our condition to get \reff{et1} is equivalent to
\[
\F(0)<p_{c}:=\inf\{p>0\,;\,\theta(p)=1\}\, , 
\]
where $\theta(p)$ is the probability that bond percolation on $\calD$ occurs with density $p$, since it is conjectured that $p_c+p_c^*=1$ (duality) for many planar graphs. In fact, by combining Corollary \ref{col1} with \reff{e6} we have:
 
\begin{corollary}\label{e1cr}
\[
1\leq p_c+p_c^*\, .
\]
\end{corollary} 

\proofof{Corollary \ref{col1}} To see this assume we have a first-passage percolation model on $\calD$ with 
\begin{equation}\label{case1}
\P(\tau_\be =0)=1-\P(\tau_\be =1)=\F(0)=1-p>p^*_c\,. 
\end{equation}
Then $\P$-a.s. there exists an infinite cluster $\calW\subseteq\calD$ composed by edges $\be$ with $\tau_\be=0$. Denote by $T(\b0,\calW)$ the first-passage time from $\b0$ to $\calW$. Then for all $t>T(\b0,\calW)$ we
have that $\bB_{\b0}(t)$ is an unbounded set. By \reff{e6} (since such a distribution satisfies \reff{c1} and \reff{c2}), this implies that $\mu(\F)=\mu(p)=0$ if $1-p>p_c$. On the other hand, by Corollary \ref{col1}, $\mu(p)>0$ if $1-p<1-p_c^*$, and so \reff{e1cr} must hold. \endproof

Other passage times have been considered in the literature such as
$T(\b0,\bH_n)$, where $\bH_n$ is the hyperplane consisting of points
$\bx=(x^1,x^2)$ so that $x_1=n$, and $T(\b0,\partial[-n,n]^2)$. The arguments in this article can be used to prove the analog of Theorem \ref{t1} when $T(\b0,\bn)$ is replaced by $T(\b0,\bH_n)$ or
$T(\b0,\partial[-n,n]^2)$. For site versions of FPP models the method works as well if we change the condition on $\F$ to $\F(0)<1-\bar{p}_c$, where now $\bar{p}_c$ is the critical probability for site percolation. Similarly to Corollary \ref{e1cr}, in this case one can also obtain the inequality $1/2\leq \bar{p}_c$. For more details we refer to \cite{P04}. 

\section{Renormalization}
For the moment we assume that $\F$ is Bernoulli with parameter $p$. Let $L\geq 1$ be a parameter whose value will be specified later. Let $\bz=(z^{1},z^{2})\in\Z^{2}$ and
\[
|\bz|_{\infty}:=\max_{j=1,2}\{|z^{j}|\}\,.
\]
Denote $\bC_{z}$ the circuit composed
by sites $\bz'\in\Z^2$ with $|\bz-\bz'|_{\infty}=2$. For each $\bA\subseteq\R^{2}$, we denote by $\partial \bA$ its boundary. For each $\bz\in\Z^{2}$ and
$r\in\{j/2\,:\,j\in\N\}$ consider the box 
\[
\bB_{z}^{rL}:= Lz+[-rL,rL]^{2}\,. 
\]

Divide $\bB_{\bz}^{L/2}$ into thirty-six sub-boxes with the same size and declare that $B_{\bz}^{L/2}$ is a {\bf full box} if all these thirty-six sub-boxes contain at least one point of $\calP$. Let 
\[
H^{L}_{\bz}:=\big[\bB_{\bz'}^{L/2}\mbox{ is a full box }\forall\,\bz'\in\bC_{\bz}\big]\,.
\]
Let $\calC_L$ be the set of all self-avoiding paths $\gamma=(\bv_j)_{1\leq j\leq h}$ in $\calD$, connecting $\partial\bB_{\bz}^{L/2}$ to $\partial\bB_{\bz}^{3L/2}$ and with $\bC_{\bv_j}\cap\bB_{\bz}^{3L/2}$ for all $j=2,\dots,h-1$. Let
\[
G^{L}_{z}:=\big[t(\gamma)\geq 1\,\forall\,\gamma\in\calC_L\big]\,.
\]
We say that $\bB_{\bz}^{L/2}$ is a {\bf good box} (or that $\bz$ is a good point) if 
\[
Y_\bz^L:=\I\big(H^{L}_{\bz}\cap G^{L}_{\bz}\big)=1\, ,
\]
where $\I\big(E\big)$ denotes the indicator function of the event $E$.
\begin{figure}[tb]
\begin{center}
\includegraphics[width=0.5\textwidth]{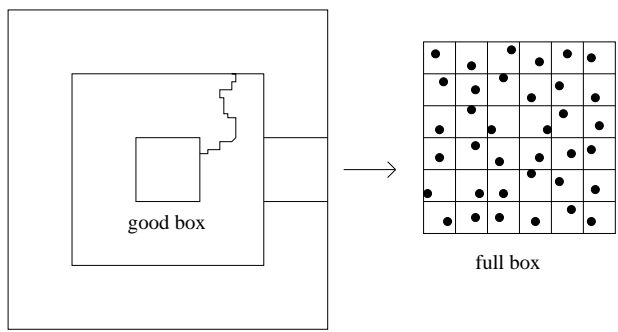}
\end{center}
\caption{Renormalization}\label{F:box}
\end{figure}

\begin{lemma}\label{l0}
 If $\P(\tau_\be=0)=1-p<1-p_{c}^*$ then 
\[
\lim_{L\to\infty}\P\big(Y_{\b0}^{L}=1\big)= 1\,.
\]
\end{lemma} 

\proofof{Lemma \ref{l0}} First notice that
\begin{equation}\label{e1re}
\P\big(Y_{\b0}^{L}=0\big)\leq \P\big((H^{L}_{\b0})^{c}\big)+\P\big((G^{L}_{\b0})^{c}\big)\,.
\end{equation}
By the definition of a two-dimensional homogeneous Poisson point process,
\begin{equation}\label{e2re}
 \lim_{L\to\infty}\P\big((H^{L}_{\b0})^{c}\big)=0\, .
\end{equation}

Now, let $X_{\be^{*}}:=\tau_{\be}$, where
$\be^{*}$ is the edge in $\calV_e$ (the Voronoi tessellation) dual to $\be$. Then $\{X_{\be^{*}}\,;\,\be^{*}\in\calV_e\}$ defines a bond percolation model on $\calV$ with law
$\P^*_p$. Consider the rectangles
\[
R^1_L:=[L/2,3L/2]\times[-3L/2,3L/2]\,,\,R^2_L:=[-3L/2,3L/2]\times[L/2,3L/2]
\]
\[
R_L^3:=[-3L/2,-L/2]\times[-3L/2,3L/2]\mbox{ and }R_L^4:=[-3L/2,3L/2]\times[-3L/2,-L/2]\,.
\]
We denote by $A_L^i$ the event $A_L$ (recall the definition of $p_c^*$) but now translate to the rectangle $R_L^i$, and by $F_{L}$ the event that an open circuit $\sigma^{*}$ in $\calV$ which surrounds $B_{\b0}^{L/2}$ and lies inside $B_{\b0}^{3L/2}$ does not exist. Thus one can easily see that 
\[
\cap_{i=1}^4 A_L^i\subseteq (F_L)^c\,. 
\]

Notice that if there exists an open circuit $\sigma^{*}$ in $\calV$ which surrounds $B_{\b0}^{L/2}$ and lies inside $B_{\b0}^{3L/2}$, then every path $\gamma$ in $\calC_L$ has an edge crossing with $\sigma^{*}$ and thus $t(\gamma)\geq 1$. Therefore,
\begin{equation}\label{e3re}
\P\big((G^{L}_{\b0})^{c}\big)\leq \P^*_p(F_{L})\leq 4\big(1-\P^*_p(A_{L})\big)\,.
\end{equation}
Since $p>p_{c}^*$, by using \reff{e1re}, \reff{e2re}, \reff{e3re} and the definition of $p_{c}^*$, we get Lemma \ref{l0}. \endproof

To obtain some sort of independence between the random variables $Y_{\bz}^L$ we shall study some geometrical aspects of Voronoi tilings. Given $\bA\subseteq\R^2$, let $\calI_\calP(\bA)$ be the sub-graph of $\calD$ composed of vertices $\bv_1$ in $\calD_v$ and edges $(\bv_2,\bv_3)$ in $\calD_e$ so that $\bC_{\bv_i}\cap\bA\neq\emptyset$ for all $i=1,2,3$.

\begin{lemma}\label{l2}
Let $L>0$ and $\bz\in\Z^2$. Assume that $\calP$ and $\calP'$ are two configurations of points so that $\calP\cap\bB_{\bz}^{5L/2}=\calP'\cap\bB_{\bz}^{5L/2}$ and that $\bB_{\bz'}^{L/2}$ is a full box with respect to $\calP$, for all $\bz'\in\bC_{\bz}$. Then $\calI_\calP(\bB_{\bz}^{3L/2})=\calI_{\calP'}(\bB_{\bz}^{3L/2})$.
\end{lemma}  

\proofof{Lemma \ref{l2}} By the definition of the Delaunay Triangulation, Lemma \ref{l2} holds if we prove that
\begin{equation}\label{e1l2}
 \bC_{\bv}(\calP)\cap\bB_{\bz}^{3L/2}\neq\emptyset\Rightarrow\bC_{\bv}(\calP)=\bC_{\bv}(\calP')\, . 
\end{equation}
To prove this we claim that 
\begin{equation}\label{e2l2}
\bC_{\bv}(\calP)\cap\bB_{\bz}^{3L/2}\neq\emptyset\Rightarrow\bC_{\bv}(\calP)\subseteq\bB_{\bz}^{2L}\,.
\end{equation}
If \reff{e2l2} does not hold then there exist $\bx_1\in\partial\bB_{\bz}^{3L/2}\cap\bC_\bv(\calP)$ and $\bx_2\in\partial\bB_{\bz}^{2L}\cap\bC_\bv(\calP)$ (by convexity of Voronoi tilings). Since every box $B^{L/2}_{\bz'}$ with $|\bz-\bz'|_\infty=2$ is a full box, there exist $\bv_1,\bv_2\in\calP$ so that
\[
 |\bv_1-\bx_1|\leq \sqrt{2}L/6\mbox{ and }|\bv_2-\bx_2|\leq \sqrt{2}L/6\,.
\]
Although, $\bx_1$ and $\bx_2$ belong to $\bC_{\bv}(\calP)$ and so
\[
 |\bv-\bx_1|\leq |\bv_1-\bx_1|\mbox{ and }|\bv-\bx_2|\leq |\bv_2-\bx_2|\, .
\]
Thus,
\[
 L/2\leq|\bx_1-\bx_2|\leq |\bx_1-\bv|+|\bx_2-\bv|\leq\sqrt{2}L/3\,,
\]
which leads to a contradiction since $\sqrt{2}/3<1/2$. By an analogous argument, one can prove that 
\begin{equation}\label{e3l2}
\bC_{\bv'}(\calP')\cap(\bB_{\bz}^{5L/2})^c\neq\emptyset\Rightarrow\bC_{\bv'}(\calP')\subseteq(\bB_{\bz}^{2L})^c\,.
\end{equation} 
Now suppose \reff{e1l2} does not hold. Without lost of generality, we may assume that there exists $\bv\in\calP$ with $\bC_{\bv}(\calP)\cap\bB_{\bz}^{3L/2}\neq\emptyset$ and $\bx\in\bC_{\bv}(\calP)$ with $\bx\not\in\bC_{\bv}(\calP')$. So  $\bx\in\bC_{\bv'}(\calP')$ for some $\bv'\in\calP'$. Although, $\calP\cap\bB_{\bz}^{5L/2}=\calP'\cap\bB_{\bz}^{5L/2}$ and then $\bv'\in(\bB_{\bz}^{5L/2})^c$, which is a contradiction with \reff{e2l2} and \reff{e3l2}. \endproof

For each $l\geq 1$, we say that the collection of random variables $\{ Y_{\bz}\,:\,\bz\in\Z^{2}\}$ is $l$-dependent if $\{Y_{\bz}\,:\,\bz\in\bA\}$ and $\{Y_{\bz}\,:\,\bz\in\bB\}$ are independent whenever  
\[
l<d_\infty(\bA,\bB):=\min\{|\bz-\bz'|_\infty\,:\,\bz\in\bA\mbox{ and }\bz'\in\bB\}\,. 
\]
Combining Lemma \ref{l2} with the translation invariance and the independence property of the Poisson point process we obtain: 
\begin{lemma}\label{l1}
For all $L>0$, $\{Y^{L}_\bz\,:\,\bz\in\Z^2\}$ is a $5$-dependent collection of identically distributed Bernoulli random variables.
\end{lemma} 

Denote $Y^L:=\{Y^{L}_{\bz}\,;\,\bz\in\Z^{2}\}$ and let $M_{m}(Y^L)$ be the maximum number of pairwise disjoint good  circuits in $\Z^{2}$, surrounding the origin and lying inside the box $[-m,m]^2$.

\begin{lemma}\label{l3}
If $\F(0)<1-p_c^*$ then there exists $L_0>0$ and $c_j=c_j(L_0)>0$ such that 
\[
\P\big(M_{m}(Y^{L_0})\leq c_1 m\big) \leq \exp(-c_2 m)\, .
\]
\end{lemma}

\proofof{Lemma \ref{l3}} Combining Lemmas \ref{l0} and \ref{l1} with and Theorem 0.0 of Ligget, Schonman and Stacey \cite{LSS97}, one gets that $Y^{L}$ is dominated from below by a collection $X^L:=\{X^{L}_{z}\,;\,z\in\Z^{2}\}$ of i.i.d. Bernoulli random variables with parameter $\rho(L)\to 1$ when $L\to\infty$. But for $\rho_L$ sufficiently close to $1$, we can chose $c>0$ sufficiently small, so that the probability of the event that $M_{m}(X^L)<cm$ decays exponentially fast with $m$ (see Chapter 3 of Grimmett \cite{G99}). Together with domination, this proves Lemma \ref{l3}. \endproof

The connection between the variable $M_{m}(Y^{L})$ and the first-passage time $T(\b0,\bn)$ is summarize by the following:

\begin{lemma}\label{ecr1}
\[
\frac{M^{L}_{nL^{-1}}}{6}\leq T(\b0,\bn)\,.
\] 
\end{lemma}

\proofof{Lemma \ref{ecr1}} We say that $(B_{\bz_{j}}^{L/2})_{1\leq j\leq h}$ is a circuit of good boxes if $(\bz_{j})_{1\leq j\leq h}$ is a good circuit in $\Z^{2}$, and that $(B_{\bz_{j}}^{L/2})_{1\leq j\leq h}$ and $(B_{\bz'_{j}}^{L/2})_{1\leq j\leq h'}$ are $l$-distant if 
\[
d_\infty\big((\bz_{j})_{1\leq j\leq k},(\bz'_{j})_{1\leq j\leq h'}\big)>l\,.
\]
Denote $M_m^L:=M_m(Y^L)$. Notice that there exist at least $(M^{L}_{nL^{-1}}/6)$ pairwise $5$-distant circuits of good boxes surrounding the origin and lying inside $[-n,n]^2\subseteq\R^2$. Therefore, every path $\gamma$ between the origin and any point outside $[-n,n]^2$ must cross at least $(M^{L}_{nL^{-1}}/6)$ $5$-distant circuits of good boxes. We claim this yields
\begin{equation}\label{ecr}
\frac{M^{L}_{nL^{-1}}}{6}\leq t(\gamma)\, . 
\end{equation}
Indeed, assume we take two $5$-distant good boxes, say $B_{\bz_1}^{L/2}$ and $B_{\bz_2}^{L/2}$, connected by a path $\gamma$ in $\calD$. Then $\gamma$ must contain two sub-paths in $\calD$, say $\bar{\gamma}_i=(\bv^i_j)_{1\leq j\leq  h_i}$ for $i=1,2$, connecting $\partial\bB_{\bz_i}^{3L/2}$ to $\partial\bB_{\bz_i}^{5L/2}$ and with  $\bC_{\bv^i_j}\cap\bB_{\bz_i}^{3L/2}$ for all $j=2,...,h_i-1$. Since $B_{\bz_1}^{L/2}$ and $B_{\bz_2}^{L/2}$ are $5$-distant good boxes, by Lemma \ref{l2}, these sub-paths must be edge disjoint. By the definition of a good box, $t(\bar{\gamma}_1)\geq 1$ and $t(\bar{\gamma}_2)\geq 1$, which yields
\[
2\leq t(\bar{\gamma}_1)+t(\bar{\gamma}_2)\leq t(\gamma)\,.
\]
By repeating this argument inductively (on the number of good boxes which are crossed by $\gamma$) one can get \reff{ecr}. Lemma \ref{ecr1} follows directly from \reff{ecr}. \endproof

\proofof{Theorem \ref{t1}} Together with Lemma \ref{ecr1}, Lemma \ref{l3} implies Theorem \ref{t1} under \reff{case1}. For the general case, assume $\F(0)=\P(\tau_{\be}=0)<1-p_{1}$. Fix $\epsilon>0$ so that $\F(\epsilon)<1-p_{c}^*$ (we can do so since $\F$ is right-continuous). Define the auxiliary process $\tau^{\epsilon}_{\be}:=\I(\tau_{\be}>\epsilon)$ and denote by $T^{\epsilon}$ the first-passage time associated to the collection $\{\tau_{\be}^{\epsilon}\,:\,\be\in\calD_e\}$. Thus $T^{\epsilon}(\b0,\bn)\leq \epsilon^{-1}T(\b0,\bn)$. Since $\tau^{\epsilon}_{\be}$ has a Bernoulli distribution with parameter $\P\big(\tau^{\epsilon}_{\be}=0\big)=\F(\epsilon)<1-p_{c}^*$, together with the previous case this yields Theorem \ref{t1}. \endproof 

\subsection*{Acknowledgment}
This work was develop during my doctoral studies at Impa and I would like to
thank my adviser, Prof. Vladas Sidoravicius, for his dedication and
encouragement during this period. I also thank the whole administrative staff of IMPA for their assistance and CNPQ for financing my doctoral studies, without which this work would have not been possible.

\end{document}